\documentclass{birkmult}
\usepackage[latin1]{inputenc}
\usepackage[T1]{fontenc}

\usepackage{graphicx,hyperref,ifpdf,color}
\usepackage{stmaryrd,amsmath,amssymb}

\newtheorem{thm}{Theorem}%[section]

\newtheorem{prp}[thm]{Proposition}
\newtheorem{conj}[thm]{Conjecture}

\theoremstyle{definition}
\newtheorem{df}[thm]{Definition}

\theoremstyle{remark}
\newtheorem{rk}[thm]{Remark}

\numberwithin{equation}{section}

\newcommand{\CC}{\ensuremath{\mathbb C}}

\newcommand{\HH}{\ensuremath{\mathbb H}}

\newcommand{\RR}{\ensuremath{\mathbb R}}
\newcommand{\TT}{\ensuremath{\mathbb T}}

\newcommand{\ZZ}{\ensuremath{\mathbb Z}}

\newcommand{\Var}{\mathop{\mathrm{Var}}}

\newcommand{\eg}{\emph{e.g.}}
\newcommand{\ie}{\emph{i.e.}}

\newcommand{\dd}{\ensuremath{\mathrm{d}}}
\newcommand{\decr}{\mathop{\makebox[0pt][l]{\kern0.5em$\downarrow$}\bigcap}}
\newcommand{\eqd}{:=}

\newcommand{\incr}{\mathop{\makebox[0pt][l]{\kern0.5em$\uparrow$}\bigcup}}

\begin{document}

\title{Is critical 2D percolation universal?}
\author{Vincent Beffara}
\address{CNRS -- UMPA \\
         École normale supérieure de Lyon \\
         46 allée d'Italie \\
         F-69364 Lyon cedex 07 \\
         France}
\email{vbeffara@ens-lyon.fr}

 % Critical phenomena


\begin{abstract}
  The aim of these notes is to explore possible ways of extending
  Smirnov's proof of Cardy's formula for critical site-percolation on
  the triangular lattice to other cases (such as bond-percolation on the
  square lattice); the main question we address is that of the choice of
  the lattice embedding into the plane which gives rise to conformal
  invariance in the scaling limit. Even though we were not able to
  produce a complete proof, we believe that the ideas presented here go
  in the right direction.
\end{abstract}

\maketitle

\section*{Introduction}

It is a strongly supported conjectured that many discrete models of
random media, such as \eg\ percolation and the Ising model, when taken
in dimension $2$ at their critical point, exhibit conformal invariance
in the scaling limit. Indeed, the \emph{universality} principle implies
that the asymptotic behavior of a critical system after rescaling
should not depend on the specific details of the underlying lattice,
and in particular it should be invariant under rotations (at least
under suitable symmetry conditions on the underlying lattice). Since
by construction a scaling limit is also invariant under rescaling, it
is natural to expect conformal invariance, as the local behavior of a
conformal map is the composition of a rotation and a rescaling.

On the other hand, conformally invariant \emph{continuous} models have
been thoroughly studied by physicists, using tools such as conformal
field theories. In 2000, Oded Schramm (\cite{schramm:UST}) introduced
a one-parameter family of continuous bidimensional random processes
which he called $SLE$ processes, as the only possible scaling limits
in this situation, under the assumption of conformal invariance;
connections between $SLE$ and CFT are now quite well understood (see
\eg\ \cite{cardy:conformal,bauer:cftsle}).

\bigskip

However, actual convergence of discrete models to $SLE$ in the scaling
limit is known for only a few models. The case on which we focus in this
paper is that of percolation. The topic of conformal invariance for
percolation has a long history --- see~\cite{langlands:percolation} and
references therein for an in-depth discussion of it.

In the case of site-percolation on the triangular lattice, it is a
celebrated result of Smirnov (\cite{smirnov:percol}) that indeed the
limit exists and is conformally invariant. While the proof is quite
simple and extremely elegant (see section~\ref{sec:smirnov} below and
references therein), it is very specific to that particular lattice,
to the point of being almost magical; it is a very natural question
to ask how it can be generalized to other cases, and in particular to
bond-percolation on the square lattice. Universality and conformal
invariance have indeed been tested numerically for percolation in
various geometries (see\ e.g.\ \cite{langlands:universality}), and
conformal invariance (assuming the existence of the limit) is known in
the case of Voronoi percolation (see\ \cite{benjamini:voronoi}).

\bigskip

In fact, it seems that the question of convergence itself has hardly
been addressed by physicists, at least in the CFT community --- a
continuous, conformally invariant object is usually the starting
point of their work rather than its outcome. Techniques such as the
renormalization group do give reason to expect the existence of a
scaling limit and of critical exponents, but they seem to not give much
insight into the emergence of rotational invariance.

This is not surprising in itself, for a very trivial reason: Take any
discrete model for which you know that there is a conformally invariant
scaling limit, say a simple random walk on $\ZZ^2$, and deform the
underlying lattice, in a linear way, so as to change the aspect ratio
of its faces. Then the scaling limit still exists (it is the image of
the previous one by the same transformation); but obviously it is not
rotationally invariant. Since all the rescaling techniques apply exactly
the same way before and after deformation, they cannot be sufficient
to derive rotational invariance. A trace of this appears in the most
general statement of the universality hypothesis (see \eg\ \cite[section
2.4]{langlands:percolation}): To paraphrase it, given any two periodic
planar graphs, the scaling limits of critical percolation on them are
conjugated by \emph{some} linear map $g$.

\bigskip

The main question we address in these notes is the following: Given a
discrete model on a doubly periodic planar graph, how to embed this
graph into the plane so as to make the scaling limit isotropic? If the
graph has additional symmetry (as for instance in the case of the square
or triangular lattices), the embedding has to preserve this symmetry;
so a restating of the same question in the terms of the universality
hypothesis would be, absent any additional symmetry for one of the two
graphs involved, can one determine the map $g$?

The most surprising thing (to me at least) about the question, besides
the fact that it appears to actually be orthogonal to the interests
of physicists in that domain, is that its answer turns out to depend
on the model considered. In other words, there is no absolute notion
of a ``conformal embedding'' of a general graph. In the case of the
simple random walk, the answer is quite easy to obtain, though it does
not seem to have appeared in the literature in the form we present it
here; in the case of percolation, I could find no reference whatsoever,
the closest being the discussion and numerical study of \emph{striated
models} in~\cite{langlands:percolation} where, instead of looking at a
different graph, the parameter $p$ in the model is chosen to depend on
the site in $\ZZ^2$ in a periodic fashion --- which admittedly is a very
related question.

\bigskip

The paper is roughly divided into two parts. In the first one, comprised
of the firs two sections, we introduce some notation and the general
framework of the approach, and we treat the case of the simple random
walk. This is enough to prove that the correct embedding is not the same
for it as for percolation; we then argue that circle packings might give
a way to answer the question in the latter case. In the second part,
which is of a more speculative nature, we investigate Smirnov's proof
in some detail, and rephrase it in such a way that its general strategy
can be applied to general triangulations. We then describe the two main
steps of a strategy that could lead to its generalization, though we
were able to perform none of the two.

\section{Notation and setup}

\subsection{The graph}

We first define the class of triangulations of the plane we are
interested in. Let $T$ be a $3$-regular finite graph of genus $1$ (\ie,
a graph that is embeddable in the torus $\TT^2 := \RR^2/\ZZ^2$ but
not in the plane, and having only vertices of degree $3$). For ease
of notation, we assume that $T$ is equipped with a fixed embedding in
$\TT^2$, which we also denote by $T$. The dual $T^*$ of $T$ (which we
also assume to be embedded in the torus once and for all) is then a
triangulation of $\TT^2$.

Let $\hat T$ (resp.\ $\hat T^*$) be the universal cover of $T$ (resp.\
$T^*$): Then $\hat T$ and $\hat T^*$ are mutually dual, infinite,
locally finite planar graphs, on which $\ZZ^2$ acts by translation. We
are interested in natural ways of embedding $\hat T$ into the complex
plane $\CC$. Let $T_i$ (the meaning of the notation will become clear
in a minute) be the embedding obtained by pulling $T$ back using the
canonical projection from $\RR^2$ to $\TT^2$ --- we will call $T_i$ the
\emph{square embedding} of $T$.

For every $\alpha\in\CC\setminus\RR$, let $\varphi_\alpha : \CC\to\CC$
be the $\RR$-linear map defined by $\varphi_\alpha(x+iy) = x+\alpha y$
(\ie, it sends $1$ to itself and $i$ to $\alpha$) and let $T_\alpha$ be
the image of $T_i$ by $\varphi _\alpha$. For lack of a better term, we
will call $T_\alpha$ the \emph{embedding of modulus $\alpha$} of $T$ in
the complex plane.

Notice that the notation $T_\alpha$ depends on the \emph{a priori}
choice of the embedding of $T$ in the flat torus; but, up to rotation
and scaling, the set of proper embeddings of $\hat T$ obtained starting
from two different embeddings of $T$ is the same, so no generality is
lost (as far as our purpose in these notes is concerned).

\bigskip

One very useful restriction on embeddings is the following:

\begin{df}
  We say that an embedding $T_\alpha$ of $\hat T$ in the complex plane
  is \emph{balanced} if each of its vertices is the barycenter (with
  equal weights) of its neighbors; or, equivalently, if the simple
  random walk on it is a martingale.
\end{df}

\begin{prp}
  Let $T$ be a $3$-regular graph of genus $1$: Then, for every
  $\alpha\in\HH$, there is a balanced embedding of $\hat T$ in the
  complex plane with modulus $\alpha$. Moreover, this embedding is
  unique up to translations of the plane.
\end{prp}

\begin{proof}
  We only give a sketch of the proof, because expanding it to a full
  proof is both straightforward and tedious. The main remark is that any
  periodic embedding which minimizes the sum $S_2$, over a period, of
  the squared lengths of its edges is balanced: Indeed, the gradient,
  with respect to the position of a given vertex, of $S_2$ is exactly
  the difference between this point and the barycenter of its neighbors.
  (This would be true in any Euclidean space.) It is easy to use a
  compactness argument to prove the existence of such a minimizer.

  To prove uniqueness up to translation is a little trickier, but since
  it is not necessary for the rest of this paper, we allow ourselves to
  give an even sketchier argument. First, one can get rid of
  translations by assuming that a fixed vertex of $\hat T$ is put at the
  origin by the embedding; the set of all possible embeddings of modulus
  $\alpha$ is then parameterized by $2(|V(T)|-1)$ real-valued
  parameters, which are the coordinates of the locations of the other
  vertices in one period of $\hat T$. In terms of these variables, $S_2$
  is polynomial of degree $2$. It is bounded below by the squared length
  of the longest edge in the embedding, which itself is bounded below,
  up to a constant depending only on the combinatorics of the graph, by
  the square of the largest of the $2(|V(T)|-1)$ parameters; so it goes
  to infinity uniformly at infinity. This implies that its Hessian
  (which is constant) is positive definite, so $S_2$ is strictly convex
  as a function of those variables. This immediately implies the
  uniqueness of the minimizer.
\end{proof}

An essential point is that, even though our proof uses Euclidean
geometry, the fact that the embedding is balanced is a linear condition.
In particular, if the embedding $T_i$ is balanced, then so are all
the other $T_\alpha$. The corresponding \emph{a priori} embedding of
$T$ itself into the flat torus $\TT^2$ (which is also unique up to
translations) will be freely referred to as \emph{the balanced embedding
of $T$ into the torus}.

\subsection{The probabilistic model}

We   will  be   interested   in  critical   site-percolation  on   the
triangulation  $T^*_\alpha$; more  specifically, the  question  we are
interested  is the  following.  Let  $\Omega$ be  a  simply connected,
smooth domain in  the complex plane, and let $A$, $B$,  $C$ and $D$ be
four points on its boundary, in that order.  For every $\delta>0$, let
$\Omega_\delta$ be the largest  connected component (in terms of graph
connectivity) of the intersection  of $\Omega$ with $\delta T_\alpha$,
and $\Omega^*_\delta$  be its  dual graph.  $\Omega_\delta$  should be
seen  as  a  discretization   of  $\Omega$  at  scale  $\delta$.   Let
$A_\delta$, $B_\delta$,  $C_\delta$ and $D_\delta$ be  the vertices of
$\Omega_\delta$   that  are  closest   to  $A$,   $B$,  $C$   and  $D$
respectively.

The model we are most interested in is critical site-percolation on
$\Omega^*_\delta$; however, most of the following considerations
remain valid for other lattice models. Let $C_\delta(\Omega,A,B,C,D)$
be the event that there is an open crossing in $\Omega^*_\delta$,
between the intervals $A_\delta B_\delta$ and $C_\delta
D_\delta$ of its boundary. Under some symmetry conditions on
$T$, Russo-Seymour-Welsh theory ensures that at criticality,
the probability of $C_\delta(\Omega,A,B,C,D)$ is bounded away
from both $0$ and $1$ as $\delta$ goes to $0$. Its limit was
conjectured by Cardy (see~\cite{cardy:formula}) using non-rigorous
arguments from conformal field theory; actual convergence was
proved, in the case of the triangular lattice (embedded in such
a way that its faces are equilateral triangles), by Smirnov
(see~\cite{smirnov:perco,beffara:easy}). We defer the statement of
the convergence to a later time. The following definition has become
standard:

\begin{df}
  Assume that, for every choice of $(\Omega,A,B,C,D)$, the probability
  of the event $C_\delta(\Omega,A,B,C,D)$ has a limit
  $f_\alpha(\Omega,A,B,C,D)$ as $\delta\to0$ --- we will refer to this
  by saying that the model \emph{has a scaling limit}. We say that the model is
  \emph{conformally invariant in the scaling limit} if, for every
  conformal map $\Phi$ from $\Omega$ to $\Phi(\Omega)$, one has
  $$f_\alpha(\Omega,A,B,C,D) = f_\alpha(\Phi(\Omega), \Phi(A), \Phi(B),
  \Phi(C), \Phi(D)).$$
  This is equivalent to saying that $f_\alpha(\Omega,A,B,C,D)$ only
  depends on the modulus of the conformal rectangle $(\Omega,A,B,C,D)$.
\end{df}
(Notice that the extension of $\Phi$ to the boundary of $\Omega$, which
is necessary for the above definition to make sense, is ensured as soon
as $\Omega$ is assumed to be regular enough.)

\section{Periodic embeddings}

\subsection{Uniqueness of the modulus}

Given $T$, it is natural to ask whether it is possible to choose a value
for $\alpha$ which provides conformal invariance in the scaling limit.
There are two possible strategies: Either give an explicit value for
which ``a miracle occurs'' (in physical terms, for which the model is
\emph{integrable} --- this is what Smirnov did in the case of the
triangular lattice), or obtain its existence in a non-constructive way
--- which is what we are trying to do here.

A reassuring fact is that, whenever such an $\alpha$ exists, it is
essentially unique:

\begin{prp}
  \label{prp:uniquealpha}
  For every graph $T$, there are either zero or two values of $\alpha$
  such that critical site-percolation on $T_\alpha^*$ is conformally
  invariant in the scaling limit. In the latter case, the two values are
  complex conjugates of each other.
\end{prp}

\begin{proof}
  The key remark is the following: Let $\beta$ be a non-real complex
  number. Since the event $C_\delta$ is defined using purely
  combinatorial features, one can push the whole picture forward through
  $\varphi_\beta$ without changing its probability. Let $\alpha' =
  \varphi_\beta(\alpha)$: $\varphi_\beta$ then transforms $\Omega$ into
  $\varphi_\beta(\Omega)$ and the lattice $T_\alpha$ into $T_{\alpha'}$.
  So, assuming convergence on both sides, one always has
  $$f_\alpha(\Omega,A,B,C,D) = f_{\alpha'} ( \varphi_\beta(\Omega),
  \varphi_\beta(A), \varphi_\beta(B), \varphi_\beta(C), \varphi_\beta(D)
  ).$$

  In the case $\beta=-i$, $\varphi_\beta$ is simply the map $z\mapsto
  \bar z$. In that case, the modulus of the conformal rectangle
  $(\varphi_{-i}(\Omega),\bar D, \bar C, \bar B, \bar A)$ is the same as
  that of $(\Omega,A,B,C,D)$, and clearly the event $C_\delta$ is
  invariant when the order of the corners is reversed. So, conformal
  invariance for $T_\alpha$ and the previous remark implies that
  $f_{\bar \alpha} (\Omega,A,B,C,D)$ still only depends on the modulus
  of the conformal rectangle --- in other words, if critical percolation
  $T_\alpha$ is conformally invariant in the scaling limit, that is also
  the case on $T_{\bar\alpha}$.

  Now assume conformal invariance in the scaling limit for two choices
  of the modulus in the upper-half plane; these moduli can always be
  written as $\alpha$ and $\alpha'=\varphi_\beta(\alpha)$ for an
  appropriate choice of $\beta\in\HH\setminus\{i\}$. Still using the
  above remark, all that is needed to arrive to a contradiction is to
  show that $f_\alpha$ does actually depend on the modulus of the
  rectangle (\ie, that it is not constant), and that there exist two
  conformal rectangles with the same modulus and whose images by
  $\varphi_\beta$ have different moduli.

  For the former point, it is enough to prove that for every choice of
  $\rho,\rho'>0$, the probability of crossing the rectangle
  $[0,\rho]\times[0,1]$ horizontally is strictly larger than that of
  crossing $[0,\rho+\rho']\times[0,1]$. This is obvious by
  Russo-Seymour-Welsh: The event that there is a vertical dual crossing
  in $\delta T_\alpha^* \cap [\rho+\delta,\rho+\rho']\times[0,1]$ is
  independent of $C_\delta([0,\rho]\times[0,1],\rho,\rho+i,i,0)$ and its
  probability is bounded below, uniformly in $\delta<\rho'/10$, by some
  positive $\varepsilon$ depending only on $\rho$ and $\rho'$. Hence,
  still assuming that the limits all exist as $\delta\to0$,
  $$f_\alpha([0,\rho+\rho']\times[0,1], \rho+\rho', \rho+\rho'+i, i, 0)
  \leqslant (1-\varepsilon) f_\alpha([0,\rho]\times[0,1], \rho, \rho+i,
  i, 0).$$

  For the latter point, assume that $\varphi_\beta$ preserves the
  equality of moduli of conformal rectangles. Let $Q=[0,1]^2$ be the
  unit square. By symmetry, the conformal rectangles $(Q,0,1,1+i,i)$ and
  $(Q,1,1+i,i,0)$ have the same modulus; on the other hand
  $\varphi_\beta(Q)$ is a parallelogram, and by our hypothesis on
  $\varphi_\beta$ it has the same modulus in both directions. This
  easily implies that it is in fact a rhombus. If now $Q'$ is the square
  with vertices $1/2$, $1+i/2$, $1/2+i$, $i/2$, $\varphi_\beta(Q')$ is
  both a rhombus (by the same argument) and a rectangle (because its
  vertices are the midpoints of the edges of $\varphi_\beta(Q)$ which is
  a rhombus). Hence $\varphi_\beta(Q')$ is a square, and so is
  $\varphi_\beta(Q)$, and in particular $\beta=\varphi_\beta(i)=i$,
  which is in contradiction with our hypothesis. 
\end{proof}

When such a pair of moduli exists, we will denote by $\alpha_T
^{\mathrm{perc}}$ the one with positive imaginary part. The same
reasoning can be done for various models, and in each case where the
scaling limit exists and is non-trivial, there will be a pair of
moduli making it conformally invariant; we will distinguish them from
each other by using the name of the model as a superscript (so that
for instance $\alpha_T^{\mathrm{RW}}$ makes the simple random walk
conformally invariant in the scaling limit --- cf.\ below).

When an argument does not depend on the specific model (as is the case
in the next subsection), we will use the generic notation $\alpha_T$ as
a placeholder.

\subsection{Obtaining $\alpha_T$ by symmetry arguments}

It should be noted that, because the value of $\alpha_T$ (when it
exists) is uniquely defined by the combinatorics of $T$, there are cases
where additional symmetry specifies its value uniquely. Indeed, assume
$\Psi$ is a graph isomorphism of $\hat T$ which is neither a translation
nor a central symmetry; for every $\alpha$, it induces a topological
isomorphism of $T_\alpha$. Assume without loss of generality that the
origin of the plane is chosen to be one of the vertices of $T_\alpha$;
let $z_0=\Psi(0)$, $z_1=\Psi(1)$ and $z_\alpha=\Psi(\alpha)$ (notice
that both $1$ and $\alpha$ are also vertices of $T_\alpha$).

\begin{figure}[ht]
  \begin{center}
    \includegraphics[scale=0.5]{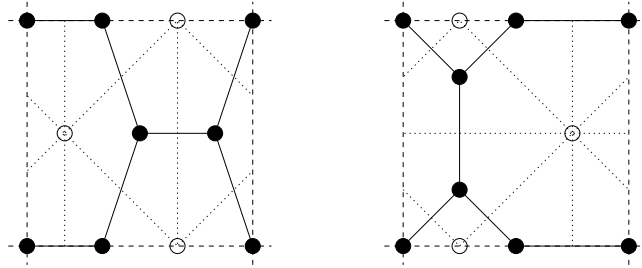}
  \end{center}
  \caption{The graphs $T_h$ (left) and $T_s$ (right), embedded into
  $\TT^2$ in a balanced way with a vertex at the origin; empty circles
  and dotted lines represent the dual graphs. Both are represented using
  their square embedding, so the triangles in $T_h$ are not
  equilateral.}
  \label{fig:hexandsquare}
\end{figure}

\begin{figure}[ht]
  \begin{center}
    \includegraphics[scale=0.5]{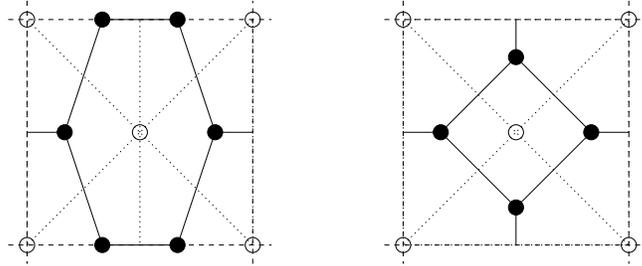}
  \end{center}
  \caption{The same graphs as in Figure~\ref{fig:hexandsquare}, with
  the origin on a vertex of the dual.}
  \label{fig:hexandsquaredual}
\end{figure}

Assume $\alpha=\alpha_T^{\mathrm{perc}}$. Because $\Psi$ is an
isomorphism, it preserves site-percolation; so, in particular, critical
site-percolation on $\Psi(T_\alpha)$ is conformally invariant in the
scaling limit. By Proposition~\ref{prp:uniquealpha}, this implies that
\begin{equation}
  \frac {z_\alpha - z_0} {z_1 - z_0} = \frac {\Psi(\alpha) - \Psi(0)}
  {\Psi(1) - \Psi(0)} \in \{ \alpha, \bar\alpha \}.
  \label{eq:psiofalpha}
\end{equation}
This condition is then enough to obtain the value of $\alpha_T$. There
are two natural examples of that (illustrated in
Figures~\ref{fig:hexandsquare} and~\ref{fig:hexandsquaredual}), which
we now describe.

\begin{itemize}
  \item Let $T_h$ be one period of the honeycomb lattice, embedded into
    $\TT^2$ in such a way that every vertex is the barycenter of its
    neighbors (we will call such an embedding \emph{balanced}); since we
    take $\TT^2$ to be a square, the coordinates of the vertices of
    $T_h$ are $(0,0)$, $(1/3,0)$, $(1/2,1/2)$ and $(5/6,1/2)$. There is
    an isomorphism $\Psi_h$ of order $3$ of $\hat T_h$, corresponding to
    rotation around $(0,0)$; on $T_\alpha$, it sends $0$ to $z_0=0$, $1$
    to $z_1=(3\alpha-1)/2$ and $\alpha$ to $z_\alpha=-(1+\alpha)/2$.
    Since $\Psi_h$ preserves orientation, Equation~\eqref{eq:psiofalpha}
    leads to $$\frac {z_\alpha - z_0} {z_1 - z_0} = \frac {1+\alpha}
    {1-3\alpha} =\alpha \quad \implies \quad \alpha = \pm i \frac
    {\sqrt3}{3},$$ in other words $\alpha_{T_h}=i\sqrt3/3$. Not
    surprisingly, this corresponds to embedding the faces of $\hat
    T_h^*$ as equilateral triangles, and those of $\hat T_h$ as regular
    hexagons.
  \item Let $T_s$ be chosen in such a way that $\hat T_s^*$ has the
    topology of the centered square lattice; if again the embedding is
    balanced, the coordinates of the vertices of $T_s$ are $(0,0)$,
    $(1/2,0)$, $(1/4,1/4)$ and $(1/4,3/4)$. There is an isomorphism
    $\Psi_s$ of order $4$ of $\hat T_s$, corresponding to a rotation
    rotation around the vertex $(1/4,0)$ of $T_s^*$. In that case
    $$\frac {z_\alpha - z_0} {z_1 - z_0} = \frac {(-3/4-\alpha/4) -
    (1/4 - \alpha/4)} {(1/4+3\alpha/4) - (1/4-\alpha/4)} = \frac
    {-1} {\alpha} = \alpha \quad \implies \quad \alpha = \pm i,$$ so
    $\alpha_{T_s}=i$. Again not surprisingly, this corresponds to the
    usual embedding of the square lattice using --- well, squares.
\end{itemize}

Of course, identifying $\alpha_T$ in those cases is a long way from a
proof of conformal invariance; but it would seem that understanding, in
the general case, what $\alpha_T^{\mathrm{perc}}$ is would be a
significant progress in our understanding of the process. 

\subsection{Embedding using random walks}

As an aside, in this subsection and the next we describe two natural
ways of embedding a doubly periodic graph into the complex plane, which
both have something to do with conformal invariance.

Let $T$ be a finite $3$-regular graph of genus $1$, embedded in $\TT^2$
in a balanced way, and let $(X_n)_{n\geqslant0}$ be a simple random
walk on it. For simplicity, assume that $(X_n)$ is irreducible as a
Markov chain. (Both $3$-regularity and irreducibility are completely
unnecessary as far as the results presented here are concerned, and the
same reasoning would work in the general case, but notation would be a
little tedious.) Since $T$ is finite, $(X_n)$ converges in distribution
to the unique invariant measure, which, because $T$ is $3$-regular, is
the uniform measure on $V(T)$; moreover the convergence is exponentially
fast.

\bigskip

Now pick $\alpha\in\HH$, and lift $(X_n)$ to a simple random walk
$(Z_n)$ on $T_\alpha$. By the balance condition on the embedding, it is
easy to check that $(Z_n)$ is a martingale; exponential decay of
correlations between its increments is enough to obtain a central limit
theorem (cf.\ for instance~\cite{jones:mcclt} and references therein).
To write the covariance matrix in a convenient form, we need some
notation. For each (oriented) edge $e$ of $T$, choose $z_1(e)$ and
$z_2(e)$ in $T_\alpha$ in such a way that they are neighbors and the
edge $(z_1(e),z_2(e))$ is a pre-image of $e$ by the natural projection
from $T_\alpha$ to $T$; let $e_\alpha:=z_2(e)-z_1(e)$ --- obviously it
does not depend on the choice of $z_1(e)$ and $z_2(e)$. Define
\begin{align*}
  \Sigma_\alpha^{xx}(T) &:= \frac 1 {|E(T)|} \sum_{e\in E(T)} (\Re
  e_\alpha)^2,\\
  \Sigma_\alpha^{yy}(T) &:= \frac 1 {|E(T)|} \sum_{e\in
  E(T)} (\Im e_\alpha)^2,\\
  \Sigma_\alpha^{xy}(T) &:= \frac 1 {|E(T)|}
  \sum_{e\in E(T)} (\Re e_\alpha)(\Im e_\alpha).
\end{align*}
It is not difficult to compute the covariance matrix of the scaling limit
of the walk:

\begin{prp}
  As $n$ goes to infinity, $n^{-1/2}Z_n$ converges in distribution to a
  Gaussian variable with covariance matrix $$\Sigma_\alpha(T) := \left[
  \begin{array}{cc} \Sigma_\alpha^{xx}(T) & \Sigma_\alpha^{xy}(T) \\
    \Sigma_\alpha^{xy}(T) & \Sigma_\alpha^{yy}(T) \end{array}\right].$$
\end{prp}

\begin{proof}
  The walk is centered by definition; the existence of a Gaussian limit
  is a direct consequence of the exponential decay of step correlations.
  All that remains to be done is to compute the covariance matrix. We
  focus on the first matrix entry, the others being similar. We have
  $$\Var(n^{-1/2}\Re Z_n) = E \left[ \frac1n \left( \sum_{k=0}^{n-1} \Re
  (Z_{k+1} - Z_k) \right)^2 \right] = \frac1n \sum_{k=0}^{n-1} E \left[
  \left( \Re (Z_{k+1} - Z_k) \right)^2 \right]$$ (the other terms
  disappear by the martingale property). We know that $Z_{k+1} - Z_k$
  converges in distribution, because the walk on $T$ converges in
  distribution; its limit is the distribution of $e_\alpha$ where $e$ is
  an edge of $T$ chosen uniformly. By Cesàro's Lemma, the expression
  above then converges to $\Sigma_\alpha^{xx}$; the computation of the
  other entries in $\Sigma_\alpha(T)$ it exactly similar.
\end{proof}

Even though the previous definition of conformal invariance in the
scaling limit does not apply directly in this case, its natural
counterpart is to ask for the scaling limit of the walk to be
rotationally invariant (\ie, to be standard two-dimensional Brownian
motion); this is equivalent to saying that the covariance matrix
$\Sigma_\alpha(T)$ is scalar, and since its entries are real, yet
another equivalent formulation is $$\left[ \Sigma_\alpha^{xx}(T) -
\Sigma_\alpha^{yy}(T) \right] + i \left[ \Sigma_\alpha^{xy}(T) \right] =
0 \quad\iff\quad \sum_{e\in E(T)} (e_\alpha)^2 = 0.$$

The last equation is a second-degree equation in $\alpha$ with
real-valued coefficients. If $\alpha\in\RR$, all the terms are
non-negative and at least one is positive, so the equation has no
solution in $\RR$; letting $\alpha$ go to $+\infty$ along the real line
leads to $|E(T)|$ positive terms, at least one of which is of order
$\alpha^2$, so the coefficient in $\alpha^2$ in the equation is not
zero. Hence the equation has exactly two solutions which are complex
conjugate of each other --- the situation is very similar to the one
in Proposition~\ref{prp:uniquealpha}. For further reference, we let
$\alpha_T^{\mathrm{RW}}$ be the one with positive imaginary part. One
advantage of this choice (besides the fact that it exists for every
doubly periodic graph) is that the value of $\alpha_T^{\mathrm{RW}}$ is
very easy to compute.

\begin{rk}
  In the more general case of a doubly periodic graph but without the
  assumptions of $3$-regularity and irreducibility (but still assuming
  that the embedding is balanced), the condition $\sum e_\alpha^2$ is
  still necessary and sufficient for the walk to be isotropic in the
  scaling limit --- and the proof is essentially the same, so we do not
  delve into more detail.
\end{rk}

\begin{rk}
  Of course, in the cases where $T$ has some additional symmetry,
  $\alpha_T^{\mathrm{RW}}$ is the same as that obtained in the previous
  subsection using symmetry \dots
\end{rk}

\begin{rk}
  One can also look at a simple random walk on the dual graph
  $T_\alpha^*$, and ask for which values of $\alpha$ this dual walk is
  isotropic in the scaling limit. As it turns out, the modulus one
  obtains this way is the same as on the initial graph, in other words
  $$\alpha_T^{\mathrm{RW}} = \alpha_{T^*}^{\mathrm{RW}}.$$ This is a
  very weak version of universality, and unfortunately there doesn't
  seem to be a purely discrete proof of it --- say, using a coupling of
  the two walks.
\end{rk}

\bigskip

There is another natural way to obtain the same condition. We are
planning on studying convergence of discrete objects to conformally
invariant limits, so it is a good idea to look for discrete-harmonic
functions on $T_\alpha$ (with respect to the natural Laplacian, which
is the same as the generator of the simple random walk on $T_\alpha$).
The condition of balanced embedding is exactly equivalent to saying that
the identity map is harmonic on $T_\alpha$; it is a linear condition, so
it does not depend on the value of $\alpha$.

The main difficulty when looking at discrete holomorphic maps is that
the product of two such maps is not holomorphic in general. But we are
interested in scaling limits, so maybe imposing that such a product is
in fact ``almost discrete holomorphic'' (in the sense that it satisfies
the Cauchy-Riemann equations up to an error term which vanishes in the
scaling limit) would be sufficient.

Whether the previous paragraph makes sense or not --- let us investigate
whether the map $\zeta : z \mapsto z^2$ is discrete-harmonic. For every
$z\in T_\alpha$, we can write $$\Delta\zeta(z) = \frac13 \sum_{z'\sim z}
(z'^2-z^2) = \frac13 \sum_{e\in E_z(T)} (z+e_\alpha)^2 - z^2 = \frac13
\sum_{e\in E_z(T)} e_\alpha^2$$ (the term in $\sum ze_\alpha$ vanishes
because the embedding is balanced). So, if $\zeta$ is discrete-harmonic,
summing the above relation over $z\in T$ gives the very same condition
$\sum e_\alpha^2=0$ as before; in other words, $\alpha_T^{\mathrm{RW}}$
is the embedding for which $z\mapsto z^2$ is \emph{discrete-harmonic
on average}.

\bigskip

As a last remark, let us investigate how strong the condition of exact
harmonicity of $\zeta$ is; so assume that $\alpha$ is chosen in such a
way that $\Delta\zeta$ is identically $0$. Let $e$ be any oriented edge
of $T$; let $e':=\tau.e$ and $e'':=\tau^2.e$ be the two other edges
sharing the same source as $e$. The balance condition on the embedding
plus harmonicity of $\zeta$ imply the following system:
\begin{equation}
  \left\{ \begin{array}{ccc}
    e_\alpha+e'_\alpha+e''_\alpha &=& 0 \\
    e_\alpha^2+(e'_\alpha)^2+(e''_\alpha)^2 &=& 0
  \end{array}\right.
  \label{eq:balanceharmonic}
\end{equation}
Up to rotation and scaling, one can always assume that $e_\alpha=1$, so
the system reduces to $e'_\alpha+e''_\alpha=-1$ and
$(e'_\alpha)^2+(e''_\alpha)^2=-1$. Squaring the first of these two
relations and substracting the second, one obtains $e'_\alpha e''_\alpha
= 1$, so $e'_\alpha$ and $e''_\alpha$ are the two solutions of the
equation $$X^2 + X + 1=0$$ which implies that $\{ e'_\alpha,e''_\alpha
\} = \{ e^{\pm 2\pi i/3} \}$. To sum it up:

\begin{prp}
  The only $3$-regular graph on which the map $\zeta:z\mapsto z^2$ is
  discrete-harmonic is the honeycomb lattice, embedded in such a way
  that its faces are regular hexagons.
\end{prp}

So, imposing $\zeta$ to be harmonic not only determines the embedding,
it also restricts $T$ to essentially one graph; but in terms of scaling
limits, the condition that $\zeta$ is harmonic on the average makes as
much sense as the exact condition.

\subsection{Embedding using circle packings}

There is another way to specify essentially unique embeddings of
triangulations, which is very strongly related to conformal geometry,
using the theory of circle packings. It is a fascinating subject in
itself and a detailed treatment would be outside of the purpose of
these notes, so the interested reader is advised to consult the book of
Stephenson~\cite{stephenson:circle} and the references therein for the
proofs of the claims in this subsection and much more.

We first give a version of a theorem of Köbe, Andreev and Thurston,
specialized to our case. It is a statement about triangulations, which
is why we actually apply it to $T^*$ instead of directly to $T$. Notice
that we \emph{do not} assume $T$ to be already embedded into the torus
$\TT^2$.

\begin{thm}[Discrete uniformization theorem
  {\cite[p.\ 51]{stephenson:circle}}]
  Let $T^*$ be a finite triangulation of the torus, and let $\hat T^*$
  be its universal cover. There exists a locally finite family
  $\smash{(\mathcal C_v)_{v\in V(\hat T^*)}}$ of disks of positive radii
  and disjoint interiors, satisfying the following compatibility
  condition: $\mathcal C_v$ and $\mathcal C_{v'}$ are tangent if, and
  only if, $v$ and $v'$ are neighbors in $\hat T^*$.
  
  Such a family is called a \emph{circle packing} associated to the
  graph $\hat T^*$. It is \emph{essentially unique}, in the following
  sense: If $(\mathcal C'_v)$ is another circle packing associated to
  $\hat T^*$, then there is a map $\varphi:\CC\to\CC$, either of the
  form $z\mapsto az+b$ or of the form $z\mapsto a\bar z+b$, such that
  for every $v\in V(\hat T^*)$, $\mathcal C'_v = \varphi(\mathcal C_v)$.
\end{thm}

\begin{rk}
  The ``existence'' part of the above theorem remains true in a much
  broader class of graphs; essentially all that is necessary is bounded
  degree and recurrence of the simple random walk on it. (One can see
  that a packing exists by completing the graph into a triangulation.)
  The ``uniqueness'' part however fails in general, as is made clear as
  soon as one tries to construct a circle packing associated to the
  square lattice \dots
\end{rk}

\bigskip

A consequence of the uniqueness part of the theorem is the following:
Let $\theta:\hat T^* \to \hat T^*$ be a translation along one of
the periods of $\hat T^*$, and let $\mathcal C'_v := \mathcal
C_{\theta(v)}$; according to the theorem, let $\varphi$ be such that
$\mathcal C'_v = \varphi(\mathcal C_v)$ for all $v$. Up to composition
of $\varphi$ by itself, one can always assume that it is of the form
$\varphi(z) = az+b$. By the assumption of local finiteness of the circle
packing, one has $|a|=1$; besides, the orbits of $\theta$ are unbounded,
so those of $\varphi$ are too, and in particular it does not have a
fixed point, which implies that $a=1$ and $b\neq0$. In other words,
$\varphi$ is a translation, \ie\ the circle packing associated to $\hat
T^*$ is itself doubly periodic.

\bigskip

\begin{figure}[ht]
  \begin{center}
    \includegraphics[scale=0.5]{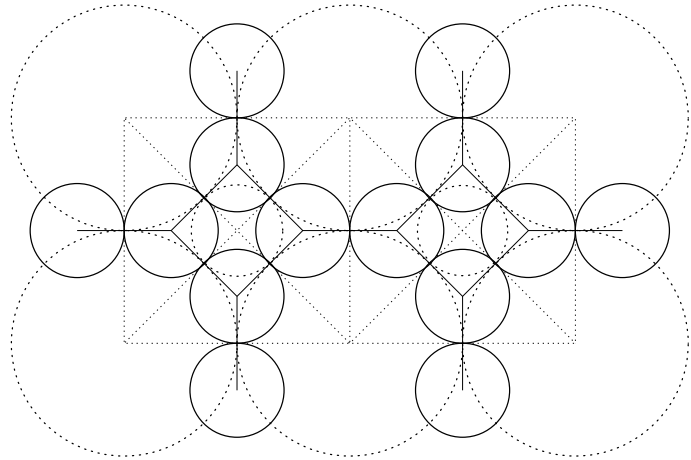}
  \end{center}
  \caption{The circle packings associated to the graph $T_s$
  (solid lines) and its dual (dotted lines).}
  \label{fig:circlepack}
\end{figure}

As soon as one is given a circle packing associated to a planar graph,
it comes with a natural embedding: Every vertex $v\in V(\hat T^*)$
will be represented by the center of $\mathcal C_v$, and if $v'$ is
a neighbor of $v$, the edge $(v,v')$ will be embedded as a segment
--- which is the union of a radius of $\mathcal C_v$ and a radius of
$\mathcal C_{v'}$, because those two disks are tangent. One can then
specify an embedding of $\hat T$ by putting each of its vertex at the
center of the disk inscribed in the corresponding triangular face of
(the embedding of) $\hat T^*$; the collection of all those inscribed
disks is in fact a circle packing associated with the graph $\hat T$
(see Figure~\ref{fig:circlepack}).

Of more interest to us is the fact that the embedding is itself doubly
periodic, by the previous remarks. Up to rotation, and scaling and
maybe complex conjugation, one can assume that the period corresponding
to the translation by $1$ (resp.\ $i$) in $\TT^2$ is equal to $1$
(resp.\ $\alpha\in\HH$). Once again, the value of the modulus $\alpha$
is uniquely determined; for further reference, we will denote it by
$\alpha_T^{\mathrm{CP}}$.

Yet again, as soon as one additional symmetry is present in $\hat T$,
the value of $\alpha_T^{\mathrm{CP}}$ is the same as that obtained
using the symmetry; this is again a direct consequence of the essential
uniqueness of the circle packing.

\subsection{``Exotic'' embeddings}

Looking closely at Smirnov's proof, one notices that essentially the
only place where the specifics of the graph are used is in the proof of
``integrability'' or exact cancellation; we will come back to this in
the next section, let us just mention that the key ingredient in the
phenomenon can be seen to be the fact that $\psi(e)$ (as introduced
earlier) is identically $0$. This is equivalent to saying that all the
triangles of the triangular lattice are equilateral.

A way to try and generalize the proof is to demand that all the faces
of $T_\alpha^*$ be equilateral triangles. Of course this cannot be done
by embedding it in the plane, even locally --- the total angle around
a vertex would be equal to $2\pi$ only if the degree of the vertex is
$6$. But one can build a $2$-dimensional manifold $M_T$ with conic
singularities by gluing together equilateral triangles according to the
combinatorics of $\hat T^*$; since the average degree of a vertex of
$T^*$ is equal to $6$, the average curvature of the manifold (defined
\eg\ as the limit of the normalized total curvature in large discs) is
$0$.

The manifold $M_T$ is not flat in general (the only case where it is
being the triangular lattice), but it is homeomorphic to the complex
plane, and one can hope to see it as a perturbation of it on which some
of the standard tools of complex analysis could have counterparts ---
the optimal being to be able to perform Smirnov's proof within it. This
is no easy task, and is probably not doable anyway.

\bigskip

To relate $M_T$ to the topic of this section, one can try to define a
module out of it. A good candidate for that is the following: Assume
that $M_T$ can be realized as a sub-manifold of $\RR^3$ (or in $\RR^d$
for $d>2$ large enough), in such a way that the (combinatorial)
translations on $\hat T$ act by global translations of the ambient
space, thus forming a \emph{periodic sub-manifold}. Then there is a
copy of $\ZZ^2$ acting on it, and the affine plane containing a given
point of $M_T$ and spanned by the directions of the two generators of
that group is at finite Hausdorff distance from it; in other words, this
realization of $M_T$ looks like a bounded perturbation of a Euclidean
plane.

One can then look at the orthogonal projections of the vertices of $\hat
T^*$ (seen as points of $M_T$) onto that plane; this creates a doubly
periodic, locally finite family of points of the Euclidean plane. It
is not always possible to form an embedding of $\hat T^*$ in the plane
from it (with disjoint edges); but it does define a value of $\alpha$ as
above.

Unfortunately, there are cases when this value of $\alpha$ is not
well-defined, in the sense that it depends on the choice of $M_T$; this
happens if the (infinite) polyhedron associated to $\hat T^*$, with
equilateral faces, is \emph{flexible}. The simplest example of this
phenomenon is to take $T$ to be two periods of the honeycomb lattice in
each direction.

\subsection{Comparing different methods of embedding}
\label{sub:comparison}

We now have at least two (forgetting about the last one) ways of giving
a conformal structure to a torus equipped with a triangulation --- which
is but another way of referring to the choice of $\alpha$. Assuming that
critical percolation does have a scaling limit, it leads to a third
choice $\alpha_T^{\mathrm{Perc}}$ of it.

\begin{figure}[ht]
  \begin{center}
    \includegraphics[scale=0.5]{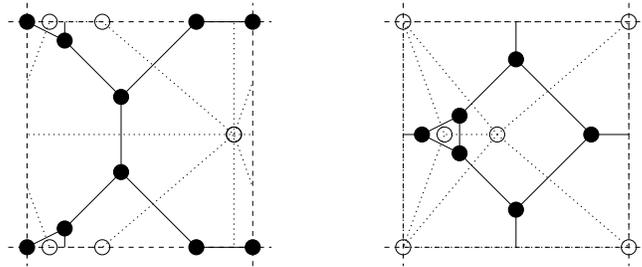}
  \end{center}
  \caption{Square (but not balanced) embeddings of $T'_s$ (solid) and
  its dual (dotted); the origin is taken as a point of $T'_s$ on the
  left, and as a point of the dual on the right, corresponding to the
  ones chosen for Figures~\ref{fig:hexandsquare} and~\ref{fig:hexandsquaredual}.}
  \label{fig:squarerefined}
\end{figure}

It would be a natural intuition that all these moduli are the same, and
correspond to a notion of \emph{conformal embedding} of a triangulation
(or a $3$-regular graph) in the plane; and they all have a claim to that
name. But this is not true in general: We detail the construction of a
counterexample. Start with the graph $T_s$ and its dual $T_s^*$; and
refine one of the ``vertical'' triangular faces of $T_s^*$ by adding a
vertex in the interior of it, connected to its three vertices. In terms
of the primal graph, this correspond to replacing one of its vertices by
a triangle --- see Figure~\ref{fig:squarerefined}. Let $T_s'$ be the
graph obtained that way; we will refer to such a splitting as a
\emph{refinement}, and to the added vertex as a \emph{new vertex}.

\begin{figure}[ht]
  \begin{center}
    \includegraphics[scale=0.5]{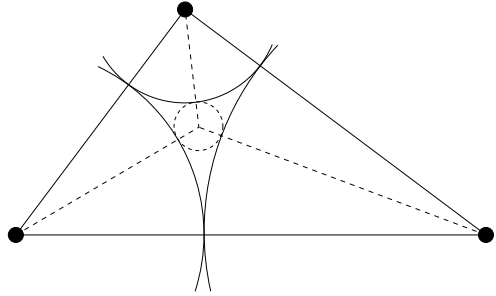}
  \end{center}
  \caption{Splitting a face of a triangulation into $3$ triangles, and
  the corresponding modification of its circle packing (added features
  represented by dashed lines)}
  \label{fig:refining}
\end{figure}

In terms of circle packings, this changes essentially nothing; the new
vertex of $(T'_s)^*$ can be realized as a new disc without modifying
the rest of the configuration (cf.\ Figure~\ref{fig:refining}).
In terms of random walks, however, adding edges will modify the
covariance matrix in the central limit theorem. The computation can be
done easily, as explained above, and one gets the following values:
$$\alpha _{T'_s} ^{\mathrm{CP}} = i = \alpha _{T_s} ^{\mathrm{CP}}
\;; \qquad \alpha _{T'_s} ^{\mathrm{RW}} = i \sqrt{\frac67} \neq i =
\alpha_{T}^{\mathrm{RW}}.$$ In this particular case, the value of
$\alpha _{(T'_s)^*} ^{\mathrm{RW}}$ is also $i\sqrt{6/7}$.

\bigskip

So, $\alpha^{\mathrm{RW}}$ and $\alpha^{\mathrm{CP}}$ are different in
general. Is $\alpha^{\mathrm{Perc}}$ (provided it exists) one of them?
An easy fact to notice is the following: Let $T^*$ be a triangulation
of the torus and let $(T')^*$ be obtained from it by splitting a
triangle into $3$ as in the construction of $T_s'$. Then, consider two
realizations of site-percolation at $p_c=1/2$ on both universal covers,
coupled in such a way that the common vertices are in the same state for
both models. In other words, start with a realization of percolation on
$\hat T^*$ and without changing site states, refine a periodic family
of triangles of it into $3$, choosing the state of each new vertex
independently of the others and of the configuration on $\hat T^*$.

If there is a chain of open vertices in $\hat T^*$, this chain is also
a chain of open vertices in the refined graph --- because all the edges
are preserved in the refinement. Conversely, starting from a chain of
open vertices in the refinement and removing each occurrence of a new
vertex on it, one obtains a chain of open vertices in $\hat T^*$; the
reason for that being that the triangle is a complete graph. Another way
of stating the same fact is to say that opening (resp.\ closing) one
of the new vertices cannot join two previously disjoint open clusters
(resp.\ split a cluster into two disjoint components); they cannot be
\emph{pivotal} for a crossing event.

Hence, the probability that a large conformal rectangle is crossed
is the same in both cases (at least if the choice of the discrete
approximation of its boundary is the same for both graphs, which in
particular implies that it contains no new vertex), and so is $f_\alpha$
for every choice of $\alpha$ (still assuming that it exists, of course).
If one is conformally invariant in the scaling limit, the other also
has to be. In short, $$\alpha _{T} ^{\mathrm{Perc}} = \alpha _{T'}
^{\mathrm{Perc}}.$$

\bigskip

Looking at circle packings instead of percolation, we get the same
identity (as was mentioned in the particular case of $T_s$), with a very
similar proof: Adding a vertex does not change anything to the rest of
the picture, and we readily obtain $$\alpha _T ^{\mathrm{CP}} = \alpha
_{T'} ^{\mathrm{CP}}.$$ This leads us to the following hope, which we
state as a conjecture even though it is much closer to being wishful
thinking:

\begin{conj}
  Let $T^*$ be a triangulation of the torus. Then, the critical
  parameter for site-percolation on its universal cover $\hat T^*$ is
  equal to $1/2$, and for every $\alpha\in\HH$, critical
  site-percolation on $\hat T^*_\alpha$ has a scaling limit. The value
  of the modulus $\alpha$ for which the model is conformally invariant
  in the scaling limit is that obtained from the circle packing
  associated to $\hat T^*$: $$\alpha_T ^{\mathrm{Perc}} = \alpha_T
  ^{\mathrm{CP}}.$$
\end{conj}

\section{Critical percolation on the triangular lattice}
\label{sec:smirnov}

For reference, and as a way of introducing our general strategy,
we give in this section a very shortened version of Smirnov's
proof of the existence and conformal invariance of the scaling
limit for critical site-percolation on the triangular lattice
$T_h$. The interested reader is advised to consult our previous
note~\cite{beffara:easy} for an ``extended shortening'', or Smirnov's
article~\cite{smirnov:percol} for the original proof; see the
book of Bollobás and Riordan~\cite{bollobas:percobook} for a more
detailed treatment. Up to cosmetic changes, we follow the notation
of~\cite{beffara:easy}.

\begin{rk}
  Up to the last paragraph of the section, we are not assuming that the
  lattice we are working with is the honeycomb lattice; our only
  assumption is that we have an \emph{a priori} bound for crossing
  probabilities of large rectangles which depends on their aspect ratio
  but not on their size (we ``assume Russo-Seymour-Welsh conditions'').
  It is not actually clear how general those are; all the standard
  proofs require at least some symmetry in the lattice in addition to
  periodicity, but it is a natural conjecture that periodicity is
  enough.
\end{rk}

\bigskip

Here and in the remainder of this paper, $\tau := e^{2\pi i/3}$ will
be the third root of unity with positive imaginary part. Let $T$ be a
finite graph of genus $1$, $T_\alpha$ an embedding of modulus $\alpha$ of
$T$ in the complex plane; let $V(T_\alpha)$ (resp.\ $E(T_\alpha)$) be
the set of vertices (resp.\ \emph{oriented} edges) of $T_\alpha$. Each
vertex $z\in V(T_\alpha)$ has three neighbours; let $E_z(T_\alpha)$ be
the set of the three oriented edges in $E(T_\alpha)$ having their source
at $z$. That set can be cyclically ordered counterclockwise; if $e\in
E_z(T_\alpha)$ is one of the three edges starting at $z$, we will denote
by $\tau.e$ (resp.\ $\tau^2.e$) the next (resp.\ second to next) edge in
the ordering.

\begin{rk}
  In the particular case of the honeycomb lattice, seeing each edge as a
  complex number (being the difference between its target and its
  source), the notation $\tau.e$ corresponds to complex multiplication
  by $e^{2\pi i/3}$ --- in other words, $\tau.e = \tau e$ as a product
  of complex numbers. That is of course not the case in general, but we
  keep the formal notation for clarity. In what follows, whenever an
  algebraic expression involves the product of a complex number by an
  edge of $T_\alpha$ or $T_\alpha^*$, as above the edge will be
  understood as the difference, as a complex number, between its target
  and its source; we will never use formal linear combinations of edges.
  The notation $\tau.e$ (with a dot) will be reserved for the
  ``topological'' rotation within $E_z(T_\alpha)$.
\end{rk}

Let again $\Omega$ be a smooth Jordan domain in the complex plane, and
let $A$, $B$, $C$ and $D$ be three points on its boundary, in that order
when following $\partial\Omega$ counterclockwise. Let $\Omega_\delta$ be
the largest connected component of $\Omega\cap\delta T_\alpha$, and let
$A_\delta$ (resp.\ $B_\delta$, $C_\delta$, $D_\delta$) be the point of
$\Omega_\delta$ that is closest to $A$ (resp.\ $B$, $C$, $D$). The main
result in Smirnov's paper~(\cite{smirnov:percol}) is the following:
\begin{thm}[Smirnov]
  \label{thm:cardy}
  In the case where $T_\alpha$ is the honeycomb lattice, embedded so as
  to make its faces regular hexagons (\ie, when $\alpha=i\sqrt3/3$),
  critical site-percolation has a conformally invariant scaling limit.
  If $\Omega$ is an equilateral triangle with vertices $A$, $B$ and $C$,
  then $$f_\alpha(\Omega,A,B,C,D) = \frac {|CD|}{|CA|}.$$
\end{thm}
Knowing this particular family of values of $f_\alpha$ is enough,
together with conformal invariance, to compute it for any conformal
rectangle. The formula obtained for a rectangle is known as
\emph{Cardy's formula}.

\bigskip

To each edge $e\in E(T_\alpha)$ corresponds its dual oriented edge $e^*
\in E(T_\alpha^*)$, oriented in such a way that the angle $(e,e^*)$
is in $(0,\pi)$. If $-e$ denotes the edge with the same endpoints as
$e$ but the reverse orientation, then we have $e^{**}=-e$. Define
$$\psi(e) := e^* + \tau (\tau.e)^* + \tau^2 (\tau^2.e)^*$$ (where as
above we interpret the edges $e^*$, $(\tau.e)^*$ and $(\tau^2.e)^*$
as complex numbers). It is easy to check that $\psi(e)=0$ if, and
only if, the face of $T_\alpha^*$ corresponding to the source of $e$
is an equilateral triangle; so, $\psi(e)$ can be seen as a measure of
the local deviation between $T_\alpha$ and the honeycomb lattice. An
identity which will be useful later is the following: \begin{equation}
\forall z\in V(T_\alpha),\quad \sum_{e\in E_z(T_\alpha)} \psi(e) = 0.
\label{eq:sumaround} \end{equation}

For every $z\in\Omega_\delta$, let $E_{A,\delta}(z)$ be the event that
there is a simple path of open vertices of $\Omega_\delta^*$, joining
two points of the boundary of the domain, which separates $z$ and $A$
from $B$ and $C$; let $H_{A} := P[E_{A,\delta}(z)]$. Define similar
events for points $B$ and $C$ by a circular permutation of the letters,
and let
\begin{align*}
 S_\delta(z) &:= H_{A,\delta}(z) + H_{B,\delta}(z)
+ H_{C,\delta}(z), \\
  H_\delta(z) &:= H_{A,\delta}(z) + \tau H_{B,\delta}(z) + \tau^2
H_{C,\delta}(z).
\end{align*}
It is a direct consequence of Russo-Seymour-Welsh estimates that these
functions are all Hölder with some universal positive exponent, with a
norm which does not depend on $\delta$, so by Ascoli's theorem they form
a relatively compact family, and as $\delta\to0$ they have subsequential
limits which are Hölder maps from $\Omega$ to $\CC$; all that is needed
is prove that only one such limit is possible.

The key argument is to show that if $h$ (resp.\ $s$) is any
subsequential limit of $(H_\delta)$ (resp.\ $(S_\delta)$) as
$\delta\to0$, then $h$ and $s$ are holomorphic; indeed, assume for
a moment that they are. Since $s$ is also real-valued, it has to be
constant, and its value is $1$ by boundary conditions (\eg\ at point
$A$). On the other hand, along the boundary arc $(A_\delta B_\delta)$
of $\partial\Omega_\delta$, $H_{C,\delta}$ is identically $0$, so the
image of the arc $(AB)$ by $h$ is contained in the segment $[1,\tau]$ of
$\CC$; and similar statements hold \emph{mutatis mutandis} for the arcs
$(BC)$ and $(CA)$. By basic index theory, this implies that $h$ is the
unique conformal map sending $\Omega$ to the (equilateral) triangle of
vertices $1$, $\tau$ and $\tau^2$, and that is enough to characterize it
and to finish the proof of Theorem~\ref{thm:cardy}.

\bigskip

So, the crux of the matter, as expected, is to prove that the map
$h$ has to be holomorphic. The most convenient way to do that is to
use Morera's theorem, which states that $h$ is indeed holomorphic on
$\Omega$ if, and only if, its integral along any closed, smooth curve
contained in $\Omega$ is equal to $0$.

Let $\gamma$ be such a curve, and let $\gamma_\delta = (z_0, z_1, \dots,
z_{L_\delta}=z_0)$ be a closed chain of vertices of $\Omega_\delta$
which approximates it within Hausdorff distance $\delta$ and has
$\mathcal O(\delta^{-1})$ points. Because the functions $H_{\delta}$ are
uniformly Hölder, it follows that $$\oint_{\gamma_\delta} H_\delta(z)
\dd z \eqd \sum _{k=0} ^{L_\delta-1} H_\delta (z_k) (z_{k+1} - z_k) \to
\oint_\gamma h(z) \dd z.$$
We want to prove that, for a suitable choice of $\alpha$, the discrete
integral on the left-hand side of that equation vanishes in the scaling
limit.

\bigskip

If $e=(z,z')$ is an oriented edge of $\Omega_\delta$, define
$P_{A,\delta}(e) := P [ E_{A,\delta}(z') \setminus E_{A,\delta}(z)]$;
define $P_{B,\delta}$ and $P_{C,\delta}$ similarly. A very clever remark
due to Smirnov, which is actually the only place in his proof where
specifics of the model (as opposed to the lattice) are used, is that one
can use color-swapping arguments to prove that, for every oriented edge,
\begin{equation}
  P_{A,\delta}(e) = P_{B,\delta}(\tau.e) = P_{C,\delta}(\tau^2.e).
  \label{eq:swapping}
\end{equation}
On the other hand, since differences of values of $H_\delta$ between
points of $\Omega_\delta$ can be computed in terms of these functions
$P_{\cdot,\delta}$, the discrete integral above can be rewritten using
them: Letting $E(\gamma_\delta)$ be the set of edges contained in the
domain surrounded by $\gamma_\delta$ and using~\eqref{eq:swapping}, one
gets
\begin{equation}
  \oint_{\gamma_\delta} H_\delta(z) \dd z = \sum_{e\in E(\gamma_\delta)}
  \psi(e) P_{A,\delta}(e) + o(1).
  \label{eq:discr}
\end{equation}
A similar computation, together with the fact that $e^* + (\tau.e)^* +
(\tau^2.e)^*$ is identically equal to $0$, leads to
\begin{equation}
  \oint_{\gamma_\delta} S_\delta(z) \dd z = o(1).
  \label{eq:discr2}
\end{equation}
We again refer the reader to~\cite{beffara:easy} for the details of this
construction.

\bigskip

Notice that it already implies that $s$ is holomorphic, hence constant
equal to $1$, independently of the value of $\alpha$; so, whether $h$
is holomorphic or not, it will send $\bar\Omega$ to the triangle of
vertices $1$, $\tau$ and $\tau^2$ anyway. In the case of the triangular
lattice embedded in the usual way, $\psi(e)$ is also identically equal
to $0$, as was mentioned above, so $h$ is itself holomorphic, and the
proof is complete.

\bigskip

The remainder of these notes is devoted to some ideas about how to
extend the general framework of the proof to more general cases; it is
not clear how close one is to a proof, but it is likely that at least
one fundamentally new idea will be required. However, we do believe that
the overall strategy which we will now describe is the right angle of
attack of the problem. Do not expect to find any formal proof in what
follows, though.

\section{Other triangulations}

\subsection{Using local shifts}

The first natural idea when trying to generalize the construction of
Smirnov is to try an apply it to more general periodic triangulations
of the plane. Indeed, in all that precedes, up to and including
Equation~\eqref{eq:discr}, nothing is specific to the regular triangular
lattice, only Russo-Seymour-Welsh conditions (and their corollary that
$p_c=1/2$) are needed. It is only at the very last step, noticing that
$\psi$ was identically equal to $0$, that the precise geometry was
needed.

The key fact that makes hope possible is the following (and it is
actually similar to one of the points we made earlier): In the
expression of the discrete integral as a sum over interior edges, each
term is the product of two contributions:
\begin{itemize}
  \item $\psi(e)$ which depends on the geometry of the embedding, and
    through that on the value of $\alpha$;
  \item $P_{A,\delta}(e)$ which is only a function of the combinatorics
    of $\Omega_\delta$.
\end{itemize}
Even though $\Omega_\delta$ as a graph does depend on the choice of
$\alpha$, one can make the following remark: Applying the transformation
$\varphi_\beta$ (for some $\beta\in\HH$) to both the domain $\Omega$ and
the lattice $\delta T_\alpha$ does not change $\Omega_\alpha$ as a
graph. In particular it does not change the value of $P_{A,\delta}(e)$.

One can then see the whole sum as a function $\beta$, say
$S_{\Omega,\delta}(\beta)$. Because $\varphi_\beta(z)$ is a real-affine
function of $\beta$, so is $S_{\Omega,\delta}$; one can then try to
solve the equation $S_{\Omega,\delta}(\beta)=0$ in $\beta$. Using the
corresponding $\varphi_\beta$, one gets a joint choice of a domain,
a lattice modulus and mesh, and a curve $\gamma$ making the discrete
contour integral vanish.

\bigskip

It the modulus thus obtained actually did not depend on $\Omega$,
$\delta$ or $\gamma$, we would be done --- call it $\alpha_T
^{\mathrm{Perc}}$ and there is only bookkeeping left to do. However we
do not even know whether it has a limit as $\delta\downarrow0$ \dots
An alternative is as follows. Because the lattice is periodic, it
makes sense to first look at the sum $\sum \psi(e) P_{A,\delta}(e)$
over one fundamental domain. If that is small, then over the copy of
the fundamental domain immediately to the right of the previous one,
the terms $\psi(e)$ are exactly the same, and one is lead to compare
$P_{A,\delta}$ for two neighboring pre-images of a given edge of $T$.

So, let $e$ be an edge of $\Omega_\delta$, and let $e+\delta$ be its
image by a translation of one period to the right. Making the dependency
on the shape of the domain explicit in the notation, one can replace the
translation of $e$ by a translation of the domain itself and the
boundary points in the opposite direction, to obtain
\begin{equation}
  P _{A,\delta} ^{\Omega} (e+\delta) = P _{A,\delta} ^{\Omega-\delta} (e).
  \label{eq:domainshift}
\end{equation}
To estimate the difference between this term and the corresponding one
in $\Omega$, one can consider coupling two realizations of percolation,
one on $\Omega_\delta$ and the other in $\Omega_\delta - \delta$, so
that they coincide on the intersection between the two.

\begin{figure}[ht]
  \begin{center}
    \input{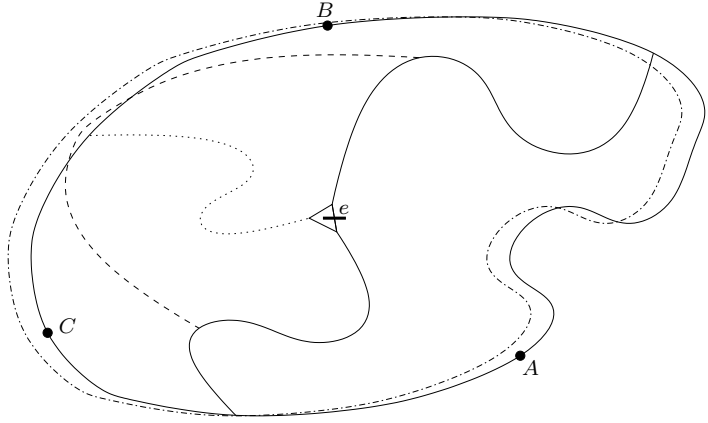}
  \end{center}
  \caption{A typical case contributing to $P_{A,\delta}
  ^{\Omega}(e) - P_{A,\delta} ^{\Omega-\delta} (e)$. The original domain
  boundary is represented by a solid line, that of the shifted domain bu
  a dashed-and-dotted line; open (resp.\ closed) arms from the source of
  $e$ are represented as solid (resp.\ dotted) lines, and the additional
  open path preventing the closed arm from connecting to the boundary of
  $\Omega_\delta-\delta$ is represented as a dashed curve.}
  \label{fig:nextarm}
\end{figure}

The event corresponding to $P_{A,\delta}^\Omega(e)$ is that there is
an open simple path separating the target of $e$ and $A$ from $B$, and
$C$, and that no open simple path separates the source of $e$ and $A$
from $B$ and $C$; this is equivalent to the existence of $3$ disjoint
paths from the $3$ vertices of the face at the source of $e$ to the
$3$ ``sides'' of the conformal triangle $(\Omega,A,B,C)$, two of them
being formed of open vertices and the third being formed of closed
vertices --- cf.\ Figure~\ref{fig:nextarm}. For this to happen in
$\Omega$ but not in $\Omega-\delta$, one of these arms needs to go up to
$\partial\Omega$ but not to $\partial(\Omega-\delta)$, and the only way
for this to be realized is for a path of the opposite color to prevent
it; this can be done in finitely many ways, Figure~\ref{fig:nextarm}
being one of them; $P_{A,\delta}^{\Omega} - P_{A,\delta}
^{\Omega-\delta}$ can then be written as the linear combination of the
probabilities of finitely many terms of that form --- half of these
actually corresponding to the reversed situation, where arms go up to
$\partial(\Omega-\delta)$ but not up to $\partial\Omega$.

In the case corresponding to Figure~\ref{fig:nextarm}, and all the
similar ones, one sees that $3$ arms connect the source of $e$ to the
boundary of $\Omega \cap (\Omega-\delta)$, and on at least one point of
that boundaries there have to be $3$ disjoint arms of diameter of order
$1$. There are $\mathcal O(\delta^{-1})$ points on the boundary, and the
probability that $3$ such arms exist from one of them is known --- at
least in the case of a polygon, which is enough for our purposes --- to
behave like $\delta^2$, see~\eg~\cite{smirnov:exps}.

Another possible reason for the non-existence of $3$ arms from the
source of $e$ to the correct portions of the boundary of $\Omega-\delta$
(say) is that one of the corresponding arms in $\Omega$ actually lands
very close to either $A$, $B$ or $C$: Preventing it from touching the
relevant part of $\partial(\Omega-\delta)$ requires only \emph{one}
additional arm from a $\delta$-neighborhood of that vertex --- \ie, a
total of $2$ arms of diameter of order $1$. The probability for that
(see~\cite{smirnov:exps} also), still in the case when $\Omega$ is a
polygon with none of $A$, $B$ or $C$ as a vertex, behaves like $\delta$.
Fortunately, there are only $3$ corners on a conformal triangle, so the
contribution of these cases is of the same order as previously

Putting everything together, one gets an estimate of the form
\begin{equation}
  P_{A,\delta}^{\Omega-\delta}(e) = P_{A,\delta}^\Omega(e) \left[ 1 +
  \mathcal O(\delta) \right].
  \label{eq:diffofpa}
\end{equation}
Coming back to our current goal, let $\mathcal E$ be the set of oriented
edges in a given period of $\Omega_\delta$, and let $\mathcal E+\delta$
be its image by the translation of vector $\delta$. Then,
\begin{align*}
  \sum_{e\in\mathcal E+\delta} \psi(e) P_{A,\delta}^\Omega (e) &= \sum
  _{e\in\mathcal E} \psi(e) P_{A,\delta}^\Omega (e) \left[ 1+\mathcal
  O(\delta) \right] \\
  &= \sum_{e\in\mathcal E} \psi(e) P_{A,\delta}^\Omega(e) + \mathcal
  O(\delta^{2+\eta})
\end{align*}
with $\eta>0$; the existence of such an $\eta$ is ensured by
Russo-Seymour-Welsh type arguments again, which ensure that, uniformly
in $e$ and $\delta$, for every edge $e$, $P_{A,\delta}^{\Omega}(e) =
\mathcal O(\delta^\eta)$.

\bigskip

Now, if that is the way the proof starts, what needs to be done is quite
clear:
\begin{itemize}
  \item Fix a period $\mathcal E$ of the graph,
  \item Choose $\alpha$ so that the previous sum, over this period, of
    $\psi(e)P_{A,\delta}^{\Omega}(e)$ is equal to $0$,
  \item Use the above estimate to give an upper bound for the same sum
    on neighboring periods;
  \item Try to somehow propagate the estimate up to the boundary.
\end{itemize}
The last part of the plan is the one that does not work directly,
because one needs of the order of $\delta^{-1}$ steps to go from
$\mathcal E$ to $\partial\Omega$, and the previous bound is not small
enough to achieve that; one would need a term of the order of $\mathcal
O(\delta^{3+\eta})$. It is however quite possible that a more careful
decomposition of the events would lead to additional cancellation,
though we were not able to perform it.

\subsection{Using incipient infinite clusters}

Another idea which might have a better chance of working out is based
on the idea of incipient infinite clusters. We are trying to ensure
that $\sum \psi(e) P_{A,\delta}(e)$ is equal to $o(\delta^2)$ over a
period for a suitable choice of $\alpha$; but for it to be exactly
equal to $0$ depends only on the ratios $P_A(e)/P_A(e')$ within the
period considered, and not on their individual values. One can then let
$\delta$ go to $0$, or equivalently let $\Omega$ increase to cover the
whole space, and look at this ratio.

\begin{prp}
  \label{prp:iic}
  There is a map $\pi:E(\hat T)\to(0,+\infty)$ such that the following
  happens. Let $e$, $e'$ be two edges of $T_\alpha$, which we identify
  with $\hat T$ for easier notation, and let $\delta=1$.  Then, as
  $\Omega$ increases to cover the whole plane, $$\frac
  {P_{A,1}^\Omega(e)} {P_{A,1}^\Omega(e')} \to \frac {\pi(e)}
  {\pi(e')},$$ uniformly in the choices of $A$, $B$ and $C$ on $\partial
  \Omega$. The map $\pi$ is periodic and does not depend on the choice
  of $\alpha$.
\end{prp}

\begin{proof}
  The argument is very similar to Kesten's proof of existence of the
  incipient infinite cluster (see~\cite{kesten:iic}); it is based on
  Russo-Seymour-Welsh estimates. It will appear in an upcoming paper
  \cite{beffara:quantitativeiic}. Notice that there is no requirement
  for $A$, $B$ and $C$ to remain separated from each other; this is
  similar to the fact that the incipient infinite cluster is also the
  limit, as $n\to\infty$, of critical percolation conditioned to the
  event that the origin is connected to the point $(n,0)$ --- which in
  turn is again a consequence of Russo-Seymour-Welsh theory. The speed
  of convergence is certainly different with and without such
  restrictions on the positions of $A$, $B$ and $C$, though.
\end{proof}

\bigskip

Seeing this Proposition, one is tempted to define $\alpha$ by solving
the equation
\begin{equation}
  \sum_{e\in\mathcal E} \psi(e) \pi(e) = 0,
  \label{eq:sumpsipi}
\end{equation}
where again the sum is taken over one period of the lattice. Indeed, all
that remains in the sum, over the same period of the lattice, of
$\psi(e) P_A(e)$ is composed of terms of a smaller order. However,
because the limit taken to define $\pi$ is uniform in the choices of
$A$, $B$ and $C$, in particular it is invariant by re-labelling of the
corners of the conformal triangle; equivalently, taking $P_B$ instead of
$P_A$ leads to the same limit. Combining this remark with
Equation~\eqref{eq:swapping}, one gets the following identities:
\begin{equation}
  \forall e \in E(\hat T), \quad \pi(e) = \pi(\tau.e) = \pi(\tau^2.e).
  \label{eq:toobadpi}
\end{equation}
In other words, $\pi(e)$ only depends on the source of $e$. For every
edge $e=(z,z')$, let $\pi(z) := \pi(e)$: If $\mathcal V$ is a period of
$V(\hat T)$, one has $$\sum_{e\in\mathcal E} \psi(e) \pi(e) =
\sum_{z\in\mathcal V} \pi(z) \sum_{z' \sim z} \psi( (z,z') ) = 0$$ by
using the remark in Equation~\eqref{eq:sumaround}.  So, the
equation~\eqref{eq:sumpsipi} is actually always true, and does not help
in finding the value of $\alpha$ \dots

\bigskip

This is actually good news, because it is the sign of emerging
cancellations in the scaling limit, which were not at first apparent;
that means that the relevant terms in~\eqref{eq:discr} are actually
smaller than they look at first sight, which in turn means that making
the leading term equal to $0$ by the correct choice of $\alpha$ leads to
even smaller terms.

Whether the overall strategy can be made to work actually depends on
the speed of convergence in the statement of Proposition~\ref{prp:iic}.
In the case of the triangular lattice, one can actually use $SLE$ to
give an explicit expansion of the ratio $P_A(e)/P_A(e')$ as $\Omega$
increases, at least in some cases; this is the subject of an upcoming
paper \cite{beffara:quantitativeiic}.

\section{Other lattices}

\subsection{Mixed percolation}

We conclude the speculative part of these notes by some considerations
about bond-percolation on the planar square lattice. The combinatorial
construction we perform here does apply to more general cases, but the
probabilistic arguments which follow do not, so we restrict ourselves to
the case of $\ZZ^2$.

\bigskip

The general idea it to map the problem of bond-percolation on $\ZZ^2$ to
one of site-percolation on a suitable triangulation of the plane. Then,
if the arguments in the previous section can be made to work, one could
potentially prove the existence and conformal invariance of a scaling
limit of critical percolation on the square lattice.

\begin{figure}[ht]
  \begin{center}
    \includegraphics[scale=.5]{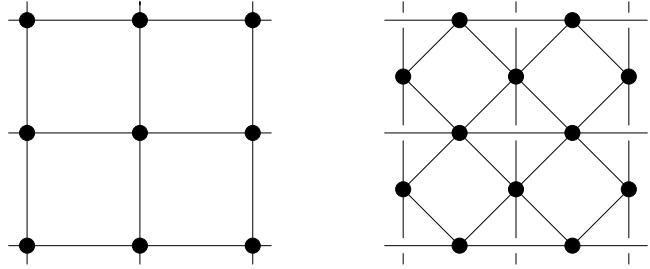}
  \end{center}
  \caption{The square lattice $\ZZ^2$ and its covering graph}
  \label{fig:covering}
\end{figure}

The key remark was already present in the book of
Kesten~\cite{kesten:book}: For any bond-percolation model on a graph,
one can construct the so-called \emph{covering graph} on which it
corresponds to site-percolation. More specifically, let $G_1$ be a
connected graph with bounded degree; as usual, let $E(G_1)$ be the set
of its edges and $V(G_1)$ be the set of its vertices. We construct a
graph $G_2$ as follows: The set $V(G_2)$ of its vertices is chosen to be
$E(G_1)$, and we put an edge between two vertices of $G_2$ if, and only
if, the corresponding edges of $G_2$ share an endpoint. Notice that even
if $G_1$ is assumed to be planar, $G_2$ does not have to be --- see
Figure~\ref{fig:covering} for the case of $\ZZ^2$.

\bigskip

The graph thus obtained from the square lattice is isomorphic to a
copy of the square lattice where every second face, in a checkerboard
disposition, is completed into a complete graph with $4$ vertices.
The next remark is the following: in terms of site-percolation, a
complete graph with $4$ vertices behaves the same way as a square with
an additional vertex at the center, which is open with probability $1$
--- with the same meaning as when we looked at refinement of triangles
in triangulations, \ie\ taking a chain of open vertices in the partially
centered square lattice and removing from it the vertices which are
face centers leads to a chain of open vertices in the covering graph of
$\ZZ^2$.

So, let again $G_s$ be the centered square lattice, as was introduced
above, and let $q\in[0,1]$; split the vertices of $G_s$ into three
classes, to defined a non-homogeneous site-percolation model, as
follows. Each vertex is either open or closed, independently of the
others, and:
\begin{itemize}
  \item The sites of $\ZZ^2$ are open with probability $p=1/2$; we will
    call them \emph{vertices of type I}, or $p$-sites for short, and
    denote by $V_1$ the set of such vertices;
  \item The vertices of coordinates $(k+1/2,l+1/2)$ with $k+l$ even are
    open with probability $q$; we will call them \emph{vertices of type
    II}, or $q$-sites for short, and denote by $V_2$ the set of such
    vertices;
  \item The vertices of coordinates $(k+1/2,l+1/2)$ with $k+l$ odd are
    open with probability $1-q$; we will call them \emph{vertices of
    type III}, or $(1-q)$-sites for short, and denote by $V_3$ the set
    of such vertices.
\end{itemize}
We will refer to that model as \emph{mixed percolation} with parameters
$p=1/2$ and $q$, and denote by $P_{1/2,q}$ the associated probability
measure. Two cases are of particular interest:
\begin{itemize}
  \item If $q=1/2$, the model is exactly critical site-percolation on
    the centered square lattice $G_s$;
  \item If $q=0$ or $q=1$ (the situation is the same in both cases up to
    a translation), from the previous remarks mixed percolation then
    corresponds to critical bond-percolation on the square lattice.
\end{itemize}
Besides, all the models obtained for $p=1/2$ are critical and satisfy
Russo-Seymour-Welsh estimates.

\subsection{Model interpolation}

We are now equipped to perform an interpolation between the models at
$q=0$ and $q=1/2$. Let $(\Omega,A,B,C,D)$ be a simply connected subset
of $\ZZ^2$ equipped with $4$ boundary points --- say, a rectangle; let
$U=U_{\Omega,A,B,C,D}$ be the event, under mixed percolation with
parameters $p=1/2$ and $q$, that there is a chain of open vertices of
$\Omega$ joining the boundary arcs $(AB)$ and $(CD)$. To estimate the
difference between the probabilities of $U$ for the two models we are
most interested in, simply write
\begin{equation}
  P_{1/2,1/2}[U] - P_{1/2,0}[U] = \int_0^{1/2} \frac \partial {\partial
  q} P_{1/2,q}[U] \,\dd q.
  \label{eq:interpol}
\end{equation}
If percolation is indeed universal, then one would expect cancellation
to occur, hopefully for each value of $q$; the optimal statement being
of the form
\begin{equation}
  \lim_{\Omega \uparrow \ZZ^2} \; \sup_{A,B,C,D \in \partial \Omega} \;
  \sup_{q\in(0,1)} \; \frac \partial {\partial q} P_{1/2,q} [U] = 0.
  \label{eq:cancellation}
\end{equation}
The main ingredient in the estimation of the derivative in $q$ is, as one
might expect, a slight generalization of Russo's formula; to state it,
we need a definition:
\begin{df}
  Consider mixed percolation on $G_s$, and let $E$ be a cylindrical
  increasing event for it (\ie, an event which depends on the state of
  finitely many vertices). Given a realization $\omega$ of the model, we
  say that a vertex $v$ is \emph{pivotal for the event $E$} if $E$ is
  realized for the configuration $\omega^v$ where $v$ is made open, and
  not realized for the configuration $\omega_v$ where $v$ is made
  closed. We will denote by $\mathrm{Piv}(E)$ the (random) set of
  pivotal vertices for $E$.
\end{df}

\begin{prp}
  \label{prp:russo}
  With the above notation, one has $$\frac \partial {\partial q}
  P_{1/2,q} [U] = E_{1/2,q} \left[ \left| \mathrm{Piv}(U) \cap \Omega
  \cap V_2 \right| - \left| \mathrm{Piv}(U) \cap \Omega \cap V_3 \right|
  \right].$$
\end{prp} 

\begin{proof}
  The argument is the same as in the proof of the usual formula
  (in the case of homogeneous percolation); we refer the reader to the
  book of Grimmett \cite{grimmett:book}.
\end{proof}

\bigskip

\begin{figure}[ht]
  \begin{center}
    \input{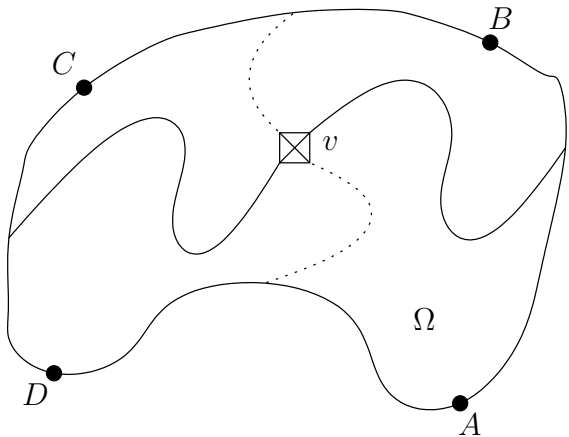}
  \end{center}
  \caption{A four-arm configuration at vertex $v$ making it pivotal for
  the event $U(\Omega,A,B,C,D)$.}
  \label{fig:fourarms}
\end{figure}

As was the case in the previous section, one can relate the event
that a given site is pivotal to the presence of disjoint arms in the
realization of the model, with appropriate color. More precisely,
a $q$-site (say) $v\in\Omega$ is pivotal if, and only if, the
following happens: $v$ is at the center of a face of $\ZZ^2$; two
opposite vertices of that face are connected respectively to the
boundary arcs $(AB)$ and $(CD)$ by disjoint chains of open vertices;
the other two vertices of the face are connected respectively to
the boundary arcs $(BC)$ and $(AD)$ by disjoint chains of closed
vertices; and none of the chains involved contains the vertex $v$.
To state the previous description more quickly, there is a $4$-arm
configuration with alternating colors at vertex $v$, and the endpoints
of the arms are appropriately located on $\partial\Omega$ --- see
Figure~\ref{fig:fourarms}.

\bigskip

The main feature of mixed percolation in the case of the centered
square lattice is the following: Starting from a configuration sampled
according to $P_{1/2,q}$ and shifting the state of all vertices by one
lattice mesh to the right, or equivalently flipping the state of all
vertices, or rotating the whole configuration by an angle of $\pi/2$
around a site of type I, one gets a configuration sampled according to
$P_{1/2,1-q}$; on the other hand, rotating the picture by $\pi/2$ around
a vertex of type II or III leaves the measure invariant.

Notice that the existence of $4$ arms of alternating colors from a given
vertex $v$ is invariant by color-swapping; the configuration in
Figure~\ref{fig:fourarms} is not though, because the arms obtained after
the color change connect the neighbors of $v$ to the wrong parts of the
boundary. Nevertheless, one can try to apply the same reasoning as in
the previous section, as follows: Let $v'$ be the vertex that is one
lattice step to the right of $v$. If $v$ is a vertex of type II, then
$v'$ is a vertex of type III, and up to boundary terms, one can pair all
the $q$-sites in $\Omega$ to corresponding $(1-q)$-sites.

To estimate the right-hand term in the statement of
Proposition~\ref{prp:russo}, let $$\Delta(v) := P[v \in \mathrm{Piv}(U)]
- P[v' \in \mathrm{Piv}(U)].$$
Our goal will be achieved if one is able to show that $\Delta(v) =
o(|\Omega|^{-1})$; or equivalently, if $\Omega_\delta$ is obtained from
a fixed continuous domain by discretization with mesh $\delta$, if one
has $$\Delta(v) = o(\delta^2).$$

\bigskip

In the case of critical site-percolation on the triangular lattice,
arguments using SLE processes give an estimate to the probability that a
vertex is pivotal, and from universality conjectures it is natural to
expect that they extend to the case of mixed percolation on $T_s$. They
involve the $4$-arm exponent of percolation, and would read (still in
the case of a fixed domain discretized at mesh $\delta$) as
$$P[v \in \mathrm{Piv}(U)] \approx \delta^{5/4}.$$

So, shifting the domain instead of the point as we did in the last
section, one would expect an estimate on $\Delta(v)$ of the order
$$\Delta(v) \approx \delta^{9/4}$$ (where the addition of $1$ in the
exponent corresponds to the presence of a $3$-arm configuration at some
point on the boundary on either the original domain or its image by the
shift). Since $9/4 > 2$, that would be enough to conclude.

However, this approach does not work directly, because of the previous
remark that the shift by one lattice step does change the measure,
replacing $q$ by $1-q$. If one is interested in the mere existence of
the $4$ arms around a vertex, combining the shift with color-flipping
is enough to cancel the effect; but the estimate one obtains that
way is of the form
\begin{equation}
  P[v \in \mathrm{Piv}(U_{\Omega,A,B,C,D})] - P[v' \in
  \mathrm{Piv}(U_{\Omega,B,C,D,A})] \approx \delta P[v \in
  \mathrm{Piv}(U_{\Omega,A,B,C,D})]
  \label{eq:almost}
\end{equation}
(and Russo-Seymour-Welsh estimates are actually enough to obtain a
formal proof of this estimate).

\bigskip

So, once again, what is missing is a way to estimate how much $P[v \in
\mathrm{Piv}(U_{\Omega,A,B,C,D})]$ depends on the location of $A$, $B$,
$C$ and $D$ along $\partial\Omega$; if the dependency is very weak, then
the estimate in Equation~\eqref{eq:almost} might actually be of the right
order of magnitude. Once again, it is likely that the way to proceed is
to use a modified version of the incipient infinite cluster conditioned
to have $4$ arms of alternating colors from the boundary, and that
the order of magnitude of $\Delta(v)$ will be related to the speed of
convergence of conditioned percolation to the incipient clusters; but we
were not able to conclude the proof that way. It would seem that this
part of the argument is easier to formalize than that of the previous
section, though, and hopefully a clever reader of these notes will be
able to do just that \dots

\bibliographystyle{amsplain}
\bibliography{Biblio}

\end{document}